\documentclass[12pt]{amsart}
\input epsf
\usepackage{amsmath}
\usepackage{amssymb}
\usepackage{amscd}
\usepackage{amsthm}
\usepackage{yhmath}
\usepackage{subfigure}
\usepackage[all]{xy}
\usepackage{color}

\numberwithin{equation}{section}

\newtheorem{theorem}[equation]{Theorem}

\newtheorem{proposition}[equation]{Proposition}

\newtheorem{lemma}[equation]{Lemma}

\newtheorem{corollary}[equation]{Corollary}

\theoremstyle{remark}
\newtheorem{remark}[equation]{Remark}

\theoremstyle{definition}
\newtheorem{definition}[equation]{Definition}

\def\XXint#1#2#3{{\setbox0=\hbox{$#1{#2#3}{\int}$} 
	\vcenter{\hbox{$#2#3$}}\kern-.5\wd0}}

\newcommand{\N}{\mathbb N}

\newcommand{\R}{\mathbb R}

\newcommand{\Z}{\mathbb Z}

\renewcommand{\S}{{\mathcal S}}

\newcommand{\diam}{\operatorname{diam}}

\newcommand{\im}{\operatorname{Im}}

\newcommand{\Isom}{\operatorname{Isom}}

\newcommand{\acts}{\curvearrowright}

\newcommand{\norm}[1]{\left\Vert#1\right\Vert}

\def\eps{\epsilon}

\def\half{{1 \over 2}}

\begin{document}

\title{Properties of isometrically homogeneous curves}

\author{Enrico Le Donne}

\date{September 4, 2011}

\begin{abstract} 
This paper is devoted to the study of isometrically homogeneous spaces
from the view point of metric geometry. Mainly we focus on those spaces that are homeomorphic to lines. 
One can reduce the study to those distances on $\R$ that are translation invariant. We study  possible values of various metric dimensions of such spaces. One of the main results is 
 the equivalence of two properties:
the first one is linear connectedness and the second one is $1$-dimensionality, with respect to  Nagata dimension. Several concrete pathological examples are provided.
\end{abstract}

\maketitle

\setcounter{secnumdepth}{2}
\setcounter{tocdepth}{2}
\tableofcontents

\section{Introduction}

An {\em isometrically homogeneous space} is a metric space $X$ with the property that the isometry group Isom$(X)$ acts on $X$ transitively.
In other words, a metric space $X$ is isometrically homogeneous if, for all $x,y\in X$, there exists an isometry $f:X\to X$  such that $f(x)=y$.

From the deep work of Gleason, Montgomery, and Zippin, it is known that, under some mild topological assumptions, an isometrically homogeneous space admits the structure of a differentiable manifold. 
Even though the topological structure of an isometrically homogeneous space is quite clear, its distance function is not well understood.
In fact, there has been no further attempt in considering a classification of isometrically homogeneous distances, in general. This paper shall make an effort in this direction.

The interest in isometrically homogeneous spaces comes from different sources.
One can consider such spaces as the analogues of topological or differentiable  homogeneous
spaces, in the metric category.
The theory of homogeneous continua has been a thriving subject for a long time.
Many structural theorems are known in low dimension, but many problems are still open. 
The curious reader might be interested in the  papers
\cite{Anderson, Rogers85, Krupski, Rogers92}.
The comprehension of differential  homogeneous spaces is a consequence of the solution of the Hilbert 5$^{th}$ problem, solved by Gleason-Montgomery-Zippin, see \cite{mz}.

Another motivation for the study of isometrically homogeneous spaces is the similarity with situations in Geometric Group Theory. 
Asymptotic cones of discrete groups are isometrically homogeneous.
Gromov hyperbolic groups act on their boundary by non-uniform biLipschitz maps transitively on a dense set.

As in the theory of homogeneous compacta, one first considers the $1$-dimensional spaces.
Namely, in this paper we mostly consider metric spaces that are topologically a line. The case when the space is topologically a circle is similar.
In fact, we shall study distances $d$ on $\R$ that are translation-invariant, i.e., 
\begin{equation}\label{translation-invariance}
d(x,y)=d(x+v,y+v), \qquad x, y, v\in \R.
\end{equation}
We call such spaces isometrically homogeneous curves. 
Indeed,
one can use Hilbert 5$^{th}$ Theory to show that such spaces are the only examples of distances on $\R$ that are isometrically homogeneous.
Actually, this fact is a special case of a more general one  about locally compact groups acting transitively on finite-dimensional locally-connected spaces.
It has been already used by Gromov \cite{Gromov-polygrowth} and Berestovski\u\i~ \cite{b2}, and it follows from results of Montgomery-Zippin \cite{mz}.
In Theorem \ref{isom_homog_curves_intro} we provide a slightly simpler proof of this result using the fact that the addition on $\R$ is the only continuous  product structure on $\R$, up to isomorphisms. 
Those distances satisfying \eqref{translation-invariance} are many and quite general ones. We shall give an equivalent description of them in Lemma \ref{characterization}.

The main focus of this paper is the study of various metric dimensions of such  isometrically homogeneous curves.
In the last part of the paper we shall present several     pathological examples. Namely, there are distances not satisfying Besicovitch Covering property, some have Nagata dimension not equal to  $1$, for some the Hausdorff dimension differs from the Assouad dimension. Before stating the formal results, in the next section we will review some notions from Metric Dimension Theory.

The deepest of our results is probably our characterization of isometrically homogeneous curves that have Nagata dimension equal to $1$.
For presenting such an equivalence, 
let us recall that a metric space $X$ is said to be {\em linearly connected} if there exists some $\lambda\geq1$ such that, for all $x,y\in X$, there exists a connected subset $S\subset X$ containing $x$ and $y$ such that $\diam(S)\leq \lambda d(x,y)$.
Linear connectedness is a standard assumption in Metric Geometry and  Geometric Group Theory, cf.  \cite{Bonk-Kleiner-05}. 
See Remark \ref{Def LC2} for an equivalent definition.
The precise notion of Nagata dimension will be recalled in the next section.
Roughly speaking, the Nagata dimension is a local-and-global analogue of the covering dimension, in the metric setting.
Our main result is the following.
\begin{theorem}\label{dimension 1 intro}
An isometrically homogeneous curve with relatively compact balls 
 has Nagata dimension $1$ if and only if it is linearly connected.
\end{theorem}
The above theorem can be made quantitative in terms of the constants of the Nagata dimension and of the linear connectedness.
Apart from some easy cases in which the Nagata dimension is $1$, in general the Nagata dimension is difficult to compute. However, one can easily check whether  a translation-invariant distance $d$ on $\R$ is linearly connected. Indeed, one just needs to verify that
$$ \sup _{t>0}\dfrac{\max\{d(0,s) \;:\; s\in(0,t)\}}{d(0,t)}<\infty.$$

\subsection{Some previous contributions and remarks}
Geodesic isometrically homogeneous spaces of finite topological dimension are well described by the work of
Berestovski\u\i~,  \cite{b2}. In fact, in the case when the distance is geodesic, then it is a Finsler-Carnot-Carath\'eodory distance with respect to a bracket-generating smooth sub-bundle.

The only other work (known to the author) on general isometrically homogeneous spaces is 
\cite{Banakh-Repovs}. The authors show that an isometrically homogeneous metric space $X$ is an Euclidean manifold if and only if $X$ is locally compact and locally
contractible. See \cite[Theorem 1.3]{Banakh-Repovs} for more general statements regarding characterization of those spaces modelled on separable or infinite-dimensional Hilbert spaces.

Examples of translation-invariant  distances on $\R$ that have Nagata dimension different from $1$ have been constructed by Higes. In fact, in
\cite{Higes} he provides examples that have infinite Nagata dimension. It would be interesting to see if there is a characterization of infinite Nagata dimensional spaces in the spirit of Theorem \ref{dimension 1 intro}.

The link between local connectedness and $1$-dimensionality was already pointed out by Freeman and Herron.
In \cite[Theorem A]{Freeman-Herron}, the authors show the following result, which resemble   Theorem \ref{dimension 1 intro}. 
Let $\gamma$ be a Jordan curve in a proper metric space.
Assume $\gamma$ is uniformly bilipschitz homogeneous.
Suppose $\gamma$ is 
not locally connected (they use the terminology `not bounded turning', which here it is the same). Then the Assouad dimension of $\gamma$ is at least $2$.
In \cite[Example 2.14]{Freeman-Herron}, they also provide an example of a proper
isometrically homogeneous   
curve that has infinite Hausdorff dimension.
The `bounded turning' property for metric Jordan curves has been also considered in \cite{Tukia-Vaisala, Meyer2011}, for characterizing (weak)-quasicircles.

	In regard to the Besicovitch Covering property, in \cite{Koranyi-Reimann, Rigot} the authors provided   examples of isometrically homogeneous spaces for which such a property does not hold.
Furthermore, in
\cite{Preiss83} Preiss claims that there are examples of metrics that are biLipschitz equivalent to the Euclidean distance on $\R$ but do not satisfy the Besicovitch Covering property.
We will give a concrete example in Theorem \ref{teoBCP0} of a distance on $\R$ that is both biLipschitz equivalent to the Euclidean one and is   translation-invariant.

\subsection{Definitions and statements of further results}
In this paper we will consider the notion of Nagata dimension. Before giving the definition, let us   recall some basic terminology. 
Two subsets $A,B$ of a metric space are \emph{$s$-separated}, for some constant $s\geq0$, if
Dist$(A, B) := \inf\{d(a,b)\;;\;a\in A,b\in B\} \geq s$.
Let  $\mathcal{B} $ be a cover of a metric space $X$. Then, for $s > 0$, the 
\emph{$s$-multiplicity} of $\mathcal{B}$ is the infimum of all $n$ such that 
every subset of $X$ with diameter $\le s$ meets at most $n$ members of the 
family  $\mathcal{B} $. 
Furthermore,
$\mathcal{B}$ is called \emph{$D$-bounded}, for some constant $D \ge 0$, 
if $\diam B := \sup\{d(x,x')\;;\;x,x'\in B\} \le D$,
for all $B \in \mathcal{B}$. 

\begin{definition}[Nagata dimension]
 Let $X$ be a metric space. The {\em Nagata dimension}, or Assouad-Nagata dimension,
 of $X$ is denoted by $\dim_N X$ and is defined as the infimum of all integers $n$ with the following property: there
exists a constant $c > 0$ such that, for all $s > 0$, $X$ admits a $cs$-bounded cover with $s$-multiplicity
at most $n + 1$.
\end{definition} 

 As shown in \cite[Proposition 2.5]{Lang-Schlichenmaier},
 the Nagata dimension  can be defined
equivalently as  the infimum of all integers $n$ with the following property: 
\begin{equation}\label{Nag Dim}
\begin{array}{l}
\text{there exists a constant $c' > 0$ such that for all $s > 0$, the metric  }\\
\text{space admits 
a $c's$-bounded covering of the form $\mathcal{B} = \bigcup_{k=0}^n \mathcal{B}_k$ }\\
\text{where 
each distinct pair of sets in  $\mathcal{B}_k$ are $s$-separated.
}
\end{array}
\end{equation}
In the next section we will 
study the  correlation between Nagata dimension  and   linear connectedness.
In fact, we shall consider the two  directions of Theorem \ref{dimension 1 intro} separately.
The difficult  implication is to show that dimension $1$ implies linear connectedness, and it will be shown in   Theorem \ref{dimension 1 one dir}. The other implication is an easy
exercise and does not require the translation invariance of the metric.
However,  the result is valid in any dimension,  under the condition that the metric is translation invariant. Recall that a metric space is called {\em proper} if it has relatively compact balls, i.e., the closed balls are compact.
In fact, in Section \ref{consequences LC} we will show the following.
\begin{theorem}\label{other dir}
Let $d$ be a proper distance on $\R^n$ that is translation invariant.
Assume that $(R^n,d)$ is linearly connected.
Then $(R^n,d)$ has Nagata dimension equal to $n$.
\end{theorem}

In section \ref{IHC}, we provide a proof for the  characterization of  isometrically homogeneous spaces homeomorphic to $\R$ as those metric spaces that are $\R$-invariant. 
The result is a corollary of results in \cite{mz} as pointed out in \cite{Gromov-polygrowth}.
Namely, we show the following.
 \begin{theorem}[Corollary of Hilbert 5$^{th}$ Theory] \label{isom_homog_curves_intro}
Let $(X,d)$ be a metric space that is topologically equivalent to $\R$.
If $\Isom(X,d)\acts X$ transitively, then 
there exists a function  $h:\R\to\R$ such that
$(X,d)$ is isometric to $(\R,d_h)$,
where
\begin{equation}
d_h(s,t):=h(|s-t|),  \qquad \forall s,t\in \R.
\end{equation}\end{theorem}

 In Lemma \ref{characterization}, for any $n\in\N$, we characterize those functions $h:\R^n\to \R$ for which the distance $d_h$ is proper, i.e., its closed balls are compact, and is $\R^n$-invariant (and of course induces the standard topology).
 
 In Section \ref{examples}, we provide four pathological examples. Before stating the next result let us 
introduce the Besicovitch Covering Property.
 %recall the definition of BCP. % and the one of Assouad dimension.
 
\begin{definition}[BCP]
 We say that a metric space $X$ satisfies
  {\em BCP} (the  {\em Besicovitch Covering Property})  if there exists an integer $N$ so that for every family of closed balls
$\mathcal A = \{A = \overline B(x_A, r_A)\}$ such that $\{x_A; A\in \mathcal A\}$ is a bounded set, one can find a
subfamily $\mathcal A'\subset  \mathcal A$ that covers the centers $x_A, A\in \mathcal A$, and so that every point in $X$
belongs to at most $N$ balls in $\mathcal A'$.
\end{definition}
See Definition \ref{defBCP} for an equivalent reformulation of such a property. See \cite{Rigot} for an historical introduction to BCP.
The first example that we provide is the following. We will deduce  it from Theorem  \ref{teoBCP1} in Remark \ref{biLipRem}.
\begin{theorem}\label{teoBCP0}
For every $L>1$ there exists a translation-invariant
distance on $\R$
that is $L$-biLipschitz equivalent to the Euclidean distance but 
 does not satisfies BCP.
 \end{theorem}

The second example has non-linearly connected balls. Thus, by Theorem \ref{dimension 1 intro}, such a space will have Nagata dimension different that the topological dimension. The construction for this result will be esplicited in Theorem \ref{teononlincon}.

\begin{theorem}\label{teononlincon0}
There exists a proper translation-invariant
distance on $\R$  such that $(\R,d)$ is not linearly connected and
has Nagata dimension $>1$. 
\end{theorem}

Other two notions of metric dimensions that we will study are the Hausdorff dimension and the  Assouad dimension. 
The Hausdorff dimension of a metric space $X$
is denoted by $\dim_H(X)$ and its
  definition     can be found in any Geometric Measure Theory, e.g., in \cite{Federer}.
We recall the Assouad dimension here.

\begin{definition}[Assouad dimension]
 The {\em Assouad dimension} of a metric space $X$ is denoted by  $\dim_{Ass}(X)$ and is defined as the infimum of all numbers $\beta > 0$ with the property that there exists some $C>1$ such that
 every set of diameter $D$ can be covered by at most
 $C\eps^{-\beta}$ sets of diameter at most $\eps D$.  
 \end{definition}

 As explained in \cite{Heinonenbook},
the Assouad dimension  can be defined
equivalently as the infimum of all numbers $\beta > 0$ with the property that every ball
of radius $r > 0$ has at most $C\eps^{-\beta}$ disjoint points of mutual distance at least $\eps r$, for
some $C> 1$ independent of the ball.

We shall provide the following two examples in Proposition \ref{Ex3}
and in in Proposition \ref{Ex4}, respectively.
\begin{theorem}\label{two examples}
There exist two proper translation-invariant
distances $d_1$, $d_2$ on $\R$  such that 
the Hausdorff dimension of $(\R,d_1)$ is infinite and
 $$\dim_{top}(\R,d_2)\neq\dim_H(\R,d_2)\neq \dim_{Ass}(\R,d_2).$$
\end{theorem}

Finally, we show an example that shows that the notion of metric dimensions that we considered are not enough to distinguish the Euclidean distance up to biLipschitz equivalence. 
Such biLipschitz equivalence can be easily checked  for translation-invariant distances, as shown in Lemma
\ref{linear rho}. Proposition \ref{dim1} will give the example.
\begin{proposition} 
There exists a  proper translation-invariant
distance $d$ on $\R$  such that 
  $(\R,d)$   has Assouad dimension $1$, Nagata dimension $1$, but it is not locally biLipschitz homeomorphic to the Euclidean line.
\end{proposition}

\section{Linear connectedness  and Nagata dimension}

This whole section is devoted to the proof of 
Theorem 
\ref{dimension 1 intro}. We begin by proving the following result.
\begin{theorem}\label{dimension 1 one dir}
Let $d$ be a proper distance on $\R$ that is translation invariant.
If $(\R,d)$ has Nagata dimension $1$, then $(\R,d)$ is linearly connected.
\end{theorem}

Theorem \ref{dimension 1 one dir}, together with Theorem \ref{other dir}
and Theorem \ref{isom_homog_curves}, will provide a proof of Theorem \ref{dimension 1 intro}.
Before showing the detailed and rigorous proof of   Theorem \ref{dimension 1 one dir}, we shall present an overview of the strategy. 

% One direction is easy and it is proved in any $\R^n$ in the next section. We show in this section the main  implication. 
\proof[Sketch of proof]
%Assume that $X:=(\R,d)$ has Nagata dimension equal to $1$ with respect to a constant $c$. 
 Assume by contradiction that $X:=(\R,d)$ is not linearly connected. Namely, we get points $0$, $x$, $y$ such that 
$0<x<y$ and $d(0,y)<<d(0,x)$.
Use $x$ and $y$ to construct a map $F$ from a ``discrete cylinder'' $C:=\Z/T\Z \times \{0,\ldots, k\}$ into $X$,
with the following property:
$${\rm Dist}( F( \Z/T\Z \times \{0\}), F(\Z/T\Z \times \{k\})) >> d(0,y),$$
$F$ is $d(0,y)$-Lipschitz, and, for some $l$, 
$$d( F(i,j), F(i+l,j) )= d(0,x), \qquad \forall i,j.$$

By assumption, $X $ has Nagata dimension equal to $1$. Namely, for some    constant $c$, we can
pick a covering of $X$ by two families $\mathcal B$ and $\mathcal W$ of sets ($\mathcal B$ is for black, $\mathcal W$ for white) such that each element in either of the families is bounded by $2 d(0,y) c$ and each two sets in the same family are $2 d(0,y)$-separated.
Now pull back the covering, i.e., an element $z$ in the cylinder $C$ is colored white if $F(z)$ belongs to an element of $\mathcal W$. The color of $z$ is black otherwise.

We shall now use the fact that any Hex game cannot end with a tie.
The fact that the black wins rephrases as the existence of a chain of `consecutive' black points in $C$ from 
$\Z/T\Z \times \{0\}$ to $\Z/T\Z \times \{k\}$. This case leads easily to a contradiction with the fact that the sets of 
$\mathcal B$ where uniformly bounded.
In case the white wins, there is a  chain of `consecutive' white points in $C$
whose associated piecewise linear curve in the `continuum' cylinder $\S^1\times [0,k]$ is closed and not null-homotopic.
Roughly speaking we get in this case the existence of two points in such a curve whose distance is approximately the distance of $0$ from $x$. Thus the image under $F$ of the chain of  white points
is contained in a single white set of big diameter.
Again, we reach a contradiction with   the boundedness of   elements in $\mathcal W$.

\subsection{Proof of   Theorem \ref{dimension 1 one dir}}
Assume that $X:=(\R,d)$ has Nagata dimension equal to $1$ with respect to a constant $c$, in the sense of Property \eqref{Nag Dim}.  Assume $c>2$.
Assume by contradiction that $X$ is not $3c$-linearly connected. Namely, using translation invariance,
we have two points   $x$, $y\in\R$ such that 
$0<x<y$ and $d(0,x)>3c d(0,y)$.
 Without loss of generality we may assume that $y=1$ and 
$d(0,y)=1$. 
We can also assume that, for some big $m\in \N$ there is $l\in\{1,2,\ldots, m-1 \}$ such that
$x=l/m$
and $$d(0,\dfrac{1}{m})\leq 1.$$
Set $R_k:=\{0, 1,\ldots, 2k(m+1) \}\times \{0,\ldots, k\}$, for $k\in\N$ to be chosen.
For a point $i\in \Z$, set $[i]$ to be the integer such that 
$[i]\in \{0,  \ldots, m \}$ and $i=[i] \mod m+1$.
Define the map $F=F_k:R_k \to X$ as
$$F(i,j):= \dfrac{[i]}{m}+j.$$
Endow $R_k$ with the $ l^1$-distance. 
We claim that $F$ is $1$-Lipschitz. Indeed, one only needs to check the images of points at distance $1$. Namely, by translation invariance,
$$d(F(i,j),(i,j+1))  = d(\dfrac{[i]}{m}+j,\dfrac{[i]}{m}+j+1)=d(0,1)=1,$$
and, if $[i]=m$,
$$d(F(i,j),(i+1,j))  = d(1+j,0+j)=d(1,0)=1,$$
if $[i]\neq m$,
$$d(F(i,j),(i+1,j))  =d(\dfrac{[i]}{m}+j,\dfrac{[i]+1}{m}+j)=d(0,\dfrac{1}{m})\leq 1.$$

We plan to chose now the value of $k$.
Notice that
$$D_0:=\{0, 1,\ldots, 2k(m+1) \}\times \{0\}$$
and
$$D_k:=\{0, 1,\ldots, 2k(m+1) \}\times \{k\}$$
have the property that
$F(D_0)\subseteq [0,1]$
and
$F(D_k)\subseteq [k,k+1]$.
Since the distance is proper, we can take $k>1$ big enough so that 
\begin{equation}\label{F D 0 F D k}
d(F(D_0),F(D_k)> 2c,
\end{equation}
where $c$ is the constant coming from the Nagata dimension, as consequence of Property \eqref{Nag Dim}.
Let  $\mathcal B$ and $\mathcal W$ be a `coloring' of $X$  by $2c$-bounded sets for which each coloring is $2$-separated.
Given a point $z\in R_k$, we say that $z$ is white if $F(z)$ belongs to an element in $\mathcal W$. We say that $z$ is black otherwise. So, if $z$ is black, $F(z)$ belongs to an element of $\mathcal B$.

We say that two points in $R_k$ are {\em neighboring} if they differ by $(1,0)$, $(1,1)$, $(0,1)$.
 Thus each point not in the `edge sides' of $R_k$ has six neighbors. Given a finite sequence of points $a_1,\ldots,a_N\in R_k$, we say that they are {\em consecutive} if $a_i$ and $a_{i-1}$ are neighboring, for $i=2,\ldots,N$.

Therefore by the property of the Hex game, see \cite{Gale}, there are two possibilities:
either there are consecutive black points 
$a_1,\ldots,a_N\in R_k$ such that $a_1\in D_1$ and $a_N\in D_k$, or there are 
 consecutive white  points 
$b_1,\ldots,b_N\in R_k$ such that $b_1\in \{0  \}\times \{0,\ldots, k\} $ and $b_N\in \{  2k(m+1) \}\times \{0,\ldots, k\}$.

Consider the first case. Let $A:=\bigcup_{i+1}^N F(a_i)$.
Since $a_i$ are consecutive, $F$ is $1$-Lipschitz, and the black family $\mathcal B$ is $2$-separated, we have that $A$ is contained in only one set of $\mathcal B$.
However, since \eqref{F D 0 F D k} and the fact that
$$A\cap F(D_0)\neq \emptyset\neq A\cap F(D_k),$$
we reach the contradiction that an element of $\mathcal B$ is not bounded by $2c$. 

Consider now the second case.
Let 
$$\pi_1: R_k=\{0, 1,\ldots, 2k(m+1) \}\times \{0,\ldots, k\}\to \{0, 1,\ldots, 2k(m+1) \},$$
$$\pi_2: R_k=\{0, 1,\ldots, 2k(m+1) \}\times \{0,\ldots, k\}\to \{0,\ldots, k\}$$
be the projections on the first and second factor of the product, respectively.
Let 
$$S_t:= \pi_1^{-1}\{t(m+1)\}, \qquad \text{ for } 0\leq t\leq 2k.$$
Clearly, the cardinality of each $S_t$ is $k+1$.
Observe that, since the $b_i$ are consecutive,
we have that
$S_t\cap\{b_1, \ldots,  b_N\}\neq \emptyset,$
for each $ 0\leq t\leq 2k.$
Since $2k>k+1$, there exist $t_1$, $t_2$,
with $t_1< t_2$, and $s_1, s_2\in  \{0,\ldots, N\}$ such that $b_{s_i}\in S_{t_i}$, $i=1,2$, and
$\pi_2(b_{s_1})=\pi_2(b_{s_2}).$

Let $\sigma$ be the piecewise linear curve on $\R^2$ determined by the points 
$b_{s_1},\ldots, b_{s_2}$.
Let $T:= m+1$. % (Actually, $T= (t_2-t_1)(m+1)$ would be good as well.)
Let $\pi:\R^2\to \R/T\Z \times \R$ be the quotient map.
Consider now the curve $\gamma:=\pi\circ \sigma$, which is closed and not null-homotopic.

To continue the proof of the theorem, we need the following algebraic-topological lemma. The proof is easy and it should be possible to find it elsewhere.
\begin{lemma}\label{rotated curve}
Let $C=\R/\Z\times \R$ be the cylinder. 
Let $\gamma$ be a closed curve in $C$ that is not null-homotopic.
For $l\in\R/\Z$, consider the new curve
$$\gamma_l: t\mapsto \gamma(t)+(l,0).$$
Then $\im (\gamma)\cap\im(\gamma_l)\neq\emptyset$ 
\end{lemma}

Using the above lemma,
there exists $p\in \im (\gamma)\cap\im(\gamma_l)$, with $l$ the value for which the initial point $x=l/m$.
Namely, we have that $p, p+(l,0)\in \im (\gamma)$.
Since the map $F$ is $(m+1,0)$-periodic, it passes to a quotient map
$$\hat F: \Z/T\Z\times \Z\to X.$$
Be aware that the point $p$ might not be a point of $\Z/T\Z\times \Z$. 
However, we can find a sequence of consecutive  points $\tilde b_{1},\ldots, \tilde b_{N}\in \im(\gamma)\cap 
\Z/T\Z\times \Z$ with $d(p,\tilde b_1)\leq1$ and $d(p+(l,0),\tilde b_N)\leq1$. 
Thus $d(\tilde b_1+(l,0),\tilde b_N)\leq 2$. 
Similarly as in the first case, the set $B:=\bigcup_{i+1}^N \hat F(\tilde b_i)$
is contained in only one set of $\mathcal W$.
However, since $c>2$,
\begin{eqnarray*}
d(\hat F(\tilde b_N),\hat F(\tilde b_1))&\geq& d(\hat F(\tilde b_1 +(l,0) ),\hat F(\tilde b_1))-2  \\
&\geq& d(\hat F(\tilde b_1) +l/m ),\hat F(\tilde b_1))-2  \\
&=& d(x,0)-2 >3c-2 > 2c.
\end{eqnarray*}
Thus we reached the contradiction that an element of $\mathcal W$ is not bounded by $2c$. 
\qed

\subsection{Consequences of linear connectedness}\label{consequences LC} 
\begin{definition}[Linear connectedness]
A metric space $X$ is said to be $\lambda$-{\em linearly connected}, for $\lambda\geq1$, if, for all $x,y\in X$, there exists a connected subset $S\subset X$ containing $x$ and $y$ such that $\diam(S)\leq \lambda d(x,y)$.
\end{definition}

\begin{remark}\label{Def LC2}
If $X$ is  $\lambda$-linearly connected, then every ball $B(p,r)\subset X$ is in the connected component of $B(p,\lambda r )$ containing $p$. On the other hand, if a metric space $X$ has the property that 
every ball $B(p,r)\subset X$ is in the connected component of $B(p,\lambda r )$ containing $p$, then $X$ is  $2\lambda$-linearly connected.
\end{remark}

In particular we point out that if $X$ is  $1$-linearly connected, then both open balls and closed balls are connected. If, moreover, $X$ is locally path connected, e.g., $\R^n$, then the open balls are path connected.

\begin{proposition}\label{1 eps l c}
Let $(X,d)$ be a   $\lambda$-linearly connected metric space.
Then there exists a metric $d'$ on $X$ that is $\lambda$-biLipschitz equivalent to $d$ and has the property that, for all $\eps>0$,  $(X,d')$ is    $(1+\eps)$-linearly connected.
\end{proposition}

\proof Define the function
\begin{equation}
d'(x,y):=\inf\{\diam(S)\;:\; \{x,y\}\subseteq S \subseteq X, S \text{ connected }\}.
\end{equation}
Such a $d'$ is a distance. Indeed, the only nontrivial property to check is the triangle inequality. For doing this, let $x,y,z\in X$ and $\eps>0.$  Take $S, T\subseteq X$ connected sets such that $ \{x,y\}\subseteq S$ and $ \{y,z\}\subseteq T$ with $\diam(S)\leq   d'(x,y)+\eps$ and  $\diam(T)\leq   d'(y,z)+\eps$
Since $S$ and $T$ intersect in $y$, their union is connected, and it contains both $x$ and $z$.
Therefore, we have
$$d'(x,z)\leq\diam(S\cup T)\leq \diam(S)+\diam(T)\leq   d'(x,y)+\eps +  d'(y,z)+\eps.$$
Hence, since $\eps$ was arbitrary, we conclude that 
$d'(x,z)\leq d'(x,y) +  d'(y,z).$

By construction, we have $d\leq d' $. By the fact that $(X,d)$ is assumed to be    $\lambda$-linearly connected,
we  $d' \leq \lambda d$. Thus  $d'$     is $\lambda$-biLipschitz equivalent to $d$.

Now we show that $d'$ is   $(1+\eps)$-linearly connected, for all $\eps>0$.
Notice that if $S$ is a connected set in $X$ then, for all $x,y\in S$, then, from the definition of $d'$,
$d'(x,y)\leq \diam(S):=\diam_d(S)$. Using also that $d\leq d'$, we conclude that
$\diam_{d'} (S)=\diam_d(S)$, for all connected sets $S$.
Hence, for all $x,y\in X$
$$d'(x,y)=\inf\{\diam_d(S)\;:\; \{x,y\}\subseteq S \subseteq X, S \text{ connected }\}$$
$$\qquad \qquad=\inf\{\diam_{d'}(S)\;:\; \{x,y\}\subseteq S \subseteq X, S \text{ connected }\}.$$
In other words,  $(X,d')$ is    $(1+\eps)$-linearly connected.
\qed

\begin{proposition}\label{1 l c}
Let $X$ be a proper metric space that is $(1+\eps)$-linearly connected, for all $\eps>0$.
Then $X$ is $1$-linearly connected. 
\end{proposition}
\proof
Let $x,y\in X$. For each $\eps>0$ we choose a  connected subset $S_\eps $ containing $x,y$ with $\diam(S_\eps)\leq (1+\eps)d(x,y)$.
Since $X$ is proper and a closure of a connected set is connected, we can assume that $S_\eps$ is compact.

By Blaschke's Theorem, there exists $\eps_n\to 0$ such that $S_{\eps_n}\to S$ with respect to the Hausdorff distance and $S$ is compact.
It is a straightforward calculation to check that $S$ is connected, contains $x$ and $y$, and $\diam(S)=d(x,y)$.
\qed

An immediate consequence of  Proposition \ref{1 l c} and the proof of Proposition \ref{1 eps l c} is the following. 
\begin{corollary}\label{corollary 1 l c}
Let $(X,d)$ be a proper, linearly connected metric space. Then there exists a biLipschitz equivalent metric $d'$ on $X$, that is proper, such that $(X,d')$ is $1$-linearly connected.
Moreover, $d'$ can be chosen such that Isom$(X,d)\subseteq$Isom$(X,d')$.
\end{corollary}

\proof[Proof of Theorem \ref{other dir}]
For Corollary \ref{corollary 1 l c} and the fact that Nagata dimension is preserved by biLipschitz map, we can also assume that $d$ is $1$-linearly connected.
In particular, each open ball is connected. Since $d$ induces the usual topology on $\R^n$, its open balls are in fact path connected.

Fix $s>0$. Consider the closure of the open ball $B:=\overline{B(0,s)}$.
Given an affine hyper-plane $\Pi$, we say that $\Pi$ is {\em tangent} to $B$ if 
$\Pi\cap\partial B\neq\emptyset$ and $\Pi\cap \stackrel{\circ}{B}=\emptyset$.
Pick $\Pi_1$ to be any hyper-plane in $\R^n$.
Let $\Pi'_1$   be an affine hyper-plane parallel to $\Pi_1$ that is tangent to $B$. 
Notice that, since $B$ is symmetric under inversion, there is only one more possibility for such a plane, namely 
 $-\Pi'_1$. 
Pick $v_1\in \Pi'_1\cap\partial B$.

By recurrence, assume we already select $v_1, \ldots, v_k$, $k<n$, linearly independent vectors in $\R^n$.
Pick a hyper-plane $\Pi_{k+1}$ (through $0$) in $\R^n$ such that
$${\rm span}\{v_1, \ldots, v_k\} \subseteq \Pi_{k}.$$
Let  $\Pi'_{k+1}$ 
be an  affine hyper-plane parallel to $\Pi_{k+1}$ tangent to $B$. 
Pick $v_{k+1}\in \Pi'_{k+1}\cap\partial B$.

Let $Q$ be  the (closed)
$n$-dimensional parallelepiped defined by the planes
$$ \Pi'_1, -\Pi'_1, \Pi'_2, -\Pi'_2, \ldots, \Pi'_n, -\Pi'_n.$$
One can construct a cover of $\R^n$ (as in the usual `brick' cover) by translations of $Q$ with the property that any other translation of $Q$ will meet at most $n+1$ elements of such a cover.
Since $B(0,s)\subseteq Q$, this cover will have $s$-multiplicity $n+1$.

The proof will be concluded if we show that $Q$ has diameter less than $cs$, for some constant $c$ depending only on $n$. For doing this, let us make use of the fact that $B$ is path connected.
For any $i=1,\ldots, n$, let $\gamma_i$ (resp. $\bar\gamma_i$) be the curve from $0$ to $v_i$ (resp. $-v_i$) inside $B$.
Let $\Gamma_i$ be the curve obtained joining together $\gamma_i$, $\bar\gamma_i$, $-\gamma_i$, $-\bar\gamma_i$, by translating them. Assume $\Gamma_i$ is parametrized by $[0,1]$.
Let $q:=-2(v_1, \ldots, v_n)$.
Consider the map
$$F:[0,1]^n\to \R^n$$
$$F(t_1, \ldots, t_n):=q+\Gamma_1(t_1)+ \ldots \Gamma_n(t_n).$$

We claim that $\diam(\im F)\leq 8 n s$ and $Q\subseteq \im F$.
Indeed,
$\diam \Gamma_i\leq 4s$ and so every point in $\im F$ has distance at most $n4s$ from $F(0)=q$.
Regarding the containment of $Q$, consider the inclusion $\iota:\partial([0,1]^n)\to[0,1]^n$.
We claim that  $\im(F\circ \iota)\cap Q=\emptyset$ and that
$$F\circ \iota :\partial([0,1]^n)\simeq S^{n+1}\to \R^n\setminus Q\simeq  S^{n+1}$$
has degree $1$, which then will give that $Q$ is covered by $\im F$.
For seeing this, we observe that $Q$ is a fundamental domain for the lattice
$$\Lambda:={\rm span}_\Z \{v_1, \ldots, v_n\}.$$
Namely,
for all $\lambda\in\Lambda\setminus\{0\}$, we have $(\lambda +\stackrel{\circ}{Q})\cap \stackrel{\circ}{Q}=\emptyset$
and
$\bigcup_{\lambda\in\Lambda}(\lambda + {Q})=\R^n$.
Moreover, the set $F    (\partial([0,1]^n))$ is in the topological annulus
$$\bigcup\{   Q+\sum_{i=1}^n a_iv_i\;:\; a_1\in\Z, |a_i|\leq 1, \sum|a_i|\neq0\}.$$
\qed

\section{Characterization of isometrically homogeneous curves}\label{IHC}
Given a function $h:\R\to\R$, define 
\begin{equation}
d_h(s,t):=h(|s-t|),
\end{equation}
for any $s,t\in\R$.
If $h$ is such that $d_h$ is a metric on $\R$, then $(\R,d_h)$ is 
an isometrically homogeneous metric space.
 Indeed, the standard translations preserve the distance $d_h$.
 
 In this section we shall  prove Theorem \ref{isom_homog_curves_intro}. Namely we show  the following characterization of dimension-one isometrically 
 homogeneous metric spaces.
 \begin{theorem}[Corollary of \cite{mz}] \label{isom_homog_curves}
Let $(X,d)$ be a metric space that is topologically equivalent to $\R$.
If $\Isom(X,d)\acts X$ transitively, then $(X,d)$ is isometric to $(\R,d_h)$
for some function $h$.
\end{theorem}
We recall again that the above theorem is not a new one. It is a special case of a general fact already pointed out in \cite{Gromov-polygrowth} and \cite{b2}.
We provide here a specific proof based on Hilbert 5$^{th}$ Theory in dimension $1$.
Here is an easy lemma needed in the proof of the theorem.
 \begin{lemma}
 In the hypotheses of the above theorem, for any $p,q\in X$, there exists a map $g\in \Isom(X,d)$ that preserves the topological orientation and is such that $g(p)=q$. Moreover, such orientation-preserving isometry is unique.
 \end{lemma}
 \proof The fact that orientation-preserving isometries act transitively is an easy connectedness argument, using the fact that the limit of orientation-preserving isometries is an orientation-preserving isometry, and that the composition of two orientation-non-preserving isometries
 is an orientation-preserving isometry.
  
  Suppose non-uniqueness, i.e., there are two such maps $g_1$ and $g_2$. Then the map $g_2^{-1}\circ g_1$ fixes $p$. We claim that the set
  $$F:=\{x\in X\;|\;(g_2^{-1}\circ g_1)(x)=x\}$$
  is all $X$, contradicting the fact that $g_1\neq g_2$.
   Trivially, $F$ is closed. Since $X$ is topologically a line, $X\setminus F$ is topologically a collection of open intervals.
   Let us identify $X$ with $\R$. So, by contradiction, if $F\neq X$, then one of the components of $X\setminus F$ is an interval $(a,b)$. So $a$ and $b$ are fixed by $f:=g_2^{-1}\circ g_1$ but not the interior of the interval.
   Now, $(a,b)\neq (-\infty,+\infty)$ since $p\in F$. We may assume $a\neq-\infty$.
   Take an auxiliary point $c\in(a,b)$ and  look at the point
   \begin{equation}
   m:=\min\{x>a\;:\; d(a,x)= d(a,c)\}
   \end{equation}
   We have $f(m)\neq m$. If $f(m)<m$, then 
   $$d(a, f(m))=d(f(a),f(m))=d(a,m)=  d(a,c),$$
  which contradicts the minimality of $m$.
    If $f(m)>m$, then    $f^{-1}(m)<m$. So
   $$d(a, f^{-1}(m))=d(f^{-1}(a),f^{-1}(m))=d(a,m)= d(a,c),$$
   which again contradicts the minimality of $m$.
   %   
%     Let us give an orientation on $X$ so we can use the left-right language.
%   Let $y$ be a right end point of one of the intervals in  $X\setminus F$, and let $r$ be the distance of the two end points of such interval.
%   Consider now all points at distance (w.r.t.~$d$) $r/2$ from $y$, (they could be more than one!).
%   If $x$ is the infimum (in the Euclidean sense) of such points, then $x$ is fixed, 
% so $x\notin X\setminus F$. Thus $F=X$.
 \qed
 
\proof[Proof of Theorem \ref{isom_homog_curves}] Set $G:=\Isom^+(X,d)$ the group of all isometries 
that preserve the topological orientation.
By the previous lemma,  $G\acts X$ transitively and, given two points $x,y\in X$, there exists a unique $g\in G$ 
with $g(x)=y$.
Therefore, fixed a point $x_0\in X$, the continuous map
\begin{eqnarray*}
\psi:&G\to&X\\
&g\mapsto&\psi(g):=g(x_0)
\end{eqnarray*}
is a bijection, and so a homeomorphism. Thus $G$ is topologically $\R$.

%Moreover, we claim that $G$ is a commutative topological group, i.e., $g\circ g'=g'\circ g$, for all $g,g'\in G$ 

Therefore, $(G,\circ)$ is a locally compact topological group homeomorphic to $\R$, (one can also easily proof that such a group is Abelian). By the solution on the Hilbert $5^{th}$ problem, $G$  is isomorphic  to $(\R,+)$, say, by a homomorphism $\phi:(\R,+)\to(G,\circ)$.
The map $$ \xymatrix{\R\ar[r]^\phi &G\ar[r]^\psi&X}$$
is a homeomorphism that we claim to be an isometry between $(\R,d_h)$ and $(X,d)$ with $h$ defined as
$$h(r):=d\left(\left(\phi(r)\right)(x_0),x_0\right).$$
Indeed,
for $s,t\in\R$, say with $s>t$,
\begin{eqnarray*} d_h(s,t)&=&h(s-t)\\
&=&d\left(\left(\phi(s-t)\right)(x_0),x_0\right)\\
&=&d\left(\left(\phi(-t+s)\right)(x_0),x_0\right)\\
%&=&d\left(\left(\phi(s)\circ\phi(t)^{-1}\right)(x_0),x_0\right)\\
&=&d\left(\left(\phi(t)^{-1}\circ\phi(s)\right)(x_0),x_0\right)\\
&=&d( \phi(s)(x_0),\phi(t)(x_0) )\\
&=&d( \psi(\phi(s)), \psi(\phi(t)) ).
\end{eqnarray*}
\qed

\section{Pathological examples}\label{examples}

\subsection{A homogeneous distance without Besicovitch property}
In this section we provide an explicit example of a distance on $\R$ that is translation invariant, biLipschitz equivalent to the Euclidean distance, and does not satisfy Besicovitch Covering property.
Such a property originates from Besicovitch's work.
Here we state the version pointed out by S\'everine Rigot.

\begin{definition}\label{defBCP}
A metric space $X$ satisfies BCP (the Besicovitch Covering Property) if there exists $N\in \N$ such that, for all bounded subsets $A\subseteq X$ and all functions $r:A\to \R_+$, there exists $\tilde A\subseteq A$ such that 
\begin{description}
\item[\ref{defBCP}.1] $$A\subseteq \bigcup_{a\in \tilde A}\overline B (a,r(a)) \qquad \text{and}$$
\item[\ref{defBCP}.2] $$\#\{a\in \tilde A\;:\; x\in  \overline B (a,r(a)) \}\leq N,\quad    \forall x\in X.$$
\end{description}
\end{definition}

Our distance will be constructed as the supremum of 
distances satisfying some properties.
For defining the distance we consider first the family of all translation-invariant  distances as follows.
Define 
$\mathcal D$ to be the collection of all distances $\rho$ on $\R$ such that $\rho(x+v,y+v)=\rho(x,y)$, for all $x,y,v\in \R$ and $\rho(0,x)\leq |x|,$ for all $x\in \R$.
Given $\rho\in \mathcal D$, we set $\rho(x):=\rho(0,x).$

\begin{definition}\label{teoBCP-def} 
For $\mathcal D$ as above, define 
$$d:=\sup\left\{\rho\in \mathcal D\;:\; \rho(\dfrac{1}{n})\leq \dfrac{1}{n+1}, \forall n\in \N\right\}.$$
\end{definition}

\begin{theorem}\label{teoBCP1}
The just-defined function $d$ is a distance satisfying the following properties
\begin{description}
\item[\ref{teoBCP1}.1]  $\half |x|\leq d(x)\leq |x|,$
\item[\ref{teoBCP1}.2]  $d$ is translation-invariant, and
\item[\ref{teoBCP1}.3]  $(\R,d)$ does not satisfies BCP.
\end{description}
\end{theorem}

The hardest property to show is the third. It is based on the fact that such a distance $d$ admits infinitely many non-connected balls.
\begin{lemma}\label{disc-ball}
There exists a sequence $r_n\to 0$ such that $\overline B (0,r_n)$ is disconnected.
\end{lemma}
Let us postpone the proof of the lemma and present the proof of Theorem \ref{teoBCP1}.
\proof[Proof of Theorem \ref{teoBCP1}]
Since $d$ is defined as a supremum of distances, then $d$ is a distance. Moreover, since such distances are  translation-invariant, then $d$ is  translation-invariant.

For $\rho(x,y):=\half |x-y|$, we have that 
$\rho\in \mathcal D$. Moreover,  for all $n\in \N$,
$\rho(\dfrac{1}{n})\leq \dfrac{1}{n+1}$, since $n+1\leq 2n$.
Thus $d\geq \rho$ and then  Property \ref{teoBCP1}.1 and \ref{teoBCP1}.2 are proved.
Regarding \ref{teoBCP1}.3, let $r_n$ be the sequence of Lemma \ref{disc-ball}.
Since $\overline B (0,r_n)$ is a disconnected subset of $\R$, there are $y'_n, y_n''>0$ such that $(y'_n,y_n'')$ is a connected component of  $\R\setminus \overline B (0,r_n)$.
Set $B_n:=\overline B (y_n''-r_n,r_n)$.
Thus $0\in B_n$, for all $n\in\N$, and $(y'_n-y_n'',0)$ is a connected component of $\R\setminus B_n$.
Up to passing to a subsequence of $r_n$, we assume that $y_{n+1}''-r_{n+1}\in(y'_n-y_n'',0)$.
In other words, the points $x_n:=y_n''-r_n$ are such that 
$0\in \overline B (x_n,r_n)$ and, if $x_n\in \overline B (x_m,r_m)$, then $n=m$.
Thus, the set $A:=\{x_n\;:\; n\in \N\}$ is counterexample for BCP. \qed

\begin{remark}
\label{biLipRem}
The result of Theorem \ref{teoBCP1} can be adapted to show Theorem \ref{teoBCP0}. Indeed, one modifies Definition \ref{teoBCP-def}  by requiring   the inequalities only for $n$ large enough in terms of $L$. 
\end{remark}
\proof[Proof of Lemma \ref{disc-ball}]
Set $$d_k:=\sup\left\{\rho\in \mathcal D\;:\; \rho(\dfrac{1}{n})\leq \dfrac{1}{n+1}, \forall n\in \N, n>k\right\}$$
and $$\delta_k:=\sup\left\{\rho\in \mathcal D\;:\; \rho(\dfrac{1}{n})\leq \dfrac{1}{n+1}, \forall n\in \N, n\leq k\right\}.$$
Thus \begin{equation}
\label{d k delta k}
d=d_k \wedge \delta_k ,\quad \forall k.\end{equation} 
Here and later, the symbol $\wedge$ denotes the minimum of two functions.

Furthermore, for all $n\in \N$, we have the implication
 \begin{equation}
\label{implication}
p\in\left(0,\dfrac{1}{n}\right)\;\text{ with }\;\delta_n(p)\leq \dfrac{1}{n+1}\; \implies\; p\leq \dfrac{1}{n+1},
\end{equation}
since, for $x\in(0,1/n)$, one has $\delta_n(x)=\min\{x, \frac{1}{n+1}-x+\frac{1}{n}\}.$

We claim that, if $p\in\left(\frac{1}{n+1},\frac{1}{n}\right)$ and $d(p)\leq \frac{1}{n+1}$, then $d_n(p)=d(p)$. Indeed, from \eqref{d k delta k}, we have $d_n(p) \wedge \delta_n(p)\leq  \frac{1}{n+1}$. 
In the case $d_n(p) \geq \delta_n(p)$,
we have $\delta_n(p)\leq  \frac{1}{n+1}$, and from \eqref{implication}, we get $p\leq  \frac{1}{n+1}$. 
But this is not the case, since $p$ is assumed to be in $\left(\frac{1}{n+1},\frac{1}{n}\right)$.
In the case $d_n(p) < \delta_n(p)$, from \eqref{d k delta k}, we have
$d(p)=d_n(p) \wedge \delta_n(p)=d_n(p)$. The claim is proved.

We now claim that 
$d_k(p)\geq \frac{k+1}{k+2}|p|$, for all $p\in \R$ and $k\in \N$.
Indeed, since the function $\frac{x}{x+1}$ is increasing, we have
$\frac{k+1}{k+2}\leq\frac{n}{n+1}$, for $k\leq n$.
So $ \frac{k+1}{k+2} \frac{1}{n}\leq\frac{1}{n+1}$ and the claim is proved.

Consequently, we can show that each ball $B:=\overline B (0,\frac{1}{n+1})$ is not connected.
Indeed, the point $1/n$ belongs to $B$, however we claim that 
$$B\cap \left(0,\dfrac{1}{n}\right) \subseteq \left(0, \dfrac{n+2}{(n+1)^2}\right).$$
Indeed, take $p\in\left(\frac{1}{n+1},\frac{1}{n}\right)$ with $d(p)\leq \frac{1}{n+1}$. 
From the claims before, we get
$$\dfrac{1}{n+1}\geq d(p)=d_n(p)\geq \dfrac{n+1}{n+2} p.$$
Thus $p\leq \frac{n+2}{(n+1)^2}$.
Since $\frac{n+2}{(n+1)^2}<\frac{1}{n}$, then such a ball $B$ is disconnected.
\qed

\subsection{A homogeneous distance without linearly connected balls}
In this section we provide an explicit example of a distance on $\R$ that is translation invariant and it is not linearly connected.
By Theorem \ref{dimension 1 intro},
such a distance will not have Nagata dimension equal to $1$.

The distance will be constructed similarly to the distance of the previous section. For this reason, we shall use the same terminology, e.g., for the set $\mathcal D$ and the notation $\rho(p)$ for $p\in \R$ and $\rho\in\mathcal D$.

\begin{definition}
Let $a_n$ be an increasing sequence such that $a_n \nearrow \infty$, $a_n/a_{n+1}$ is decreasing, and $a_{n+1}/a_n^2 \to \infty$.
Define $d$ as follows
$$d:=\sup\left\{\rho\in \mathcal D\;:\; \rho(a_{n+1})\leq a_n, \forall n\in \N\right\}.$$
\end{definition}

\begin{theorem}\label{teononlincon}
The just-defined function $d$ is a  translation-invariant distance such that
\begin{description}
\item[\ref{teononlincon}.1]  $(\R,d)$ is not linearly connected,
\item[\ref{teononlincon}.2]  the Nagata dimension of  $(\R,d)$ is not $1$. 
\end{description}
\end{theorem}

For the proof of the above theorem, we need the following auxiliary distance. Set
$$\partial_k:=\sup\left\{\rho\in \mathcal D\;:\; \rho(a_{k+1})\leq a_k\right\}.$$
Therefore we have the following properties:
\begin{equation}\label{prop 1}
d=\min_{k\in\N}\partial_k,
\end{equation}
\begin{equation}\label{prop 2}
\partial_k\geq \dfrac{a_k}{a_{k+1}}|\cdot|, \qquad\forall k\in\N, 
\end{equation}
and, setting $q_n:=(a_n+a_{n+1})/2$, we have that
\begin{equation}\label{prop 3}
\text{ on } (0,q_n),\qquad \partial_n=|\cdot|, \;\qquad\forall n\in N.
\end{equation}
\proof[Proof of Theorem \ref{teononlincon}]
It is easy to see that $d$ is a  translation-invariant distance. Also, the fact that $\dim_N(\R, d)\neq 1$ is a consequence of Theorem \ref{dimension 1 intro}.
We are just left to prove the non-linear connectedness. 
Consider the ball $B_n:=\overline B (0,a_n)$.
Obviously both $a_n$ and $a_{n+1}\in B_n$. Thus the midpoint $q_n:=(a_n+a_{n+1})/2$ is in the convex hull of $B_n$.

We claim that $d(q_n)/a_n \to \infty$, as $n\to\infty$. Indeed, from \eqref{prop 1}, 
\eqref{prop 2}, and \eqref{prop 3},
$$ %\begin{eqnarray*}
d(q_n) = \min_{k\in\N}\partial_k(q_n) %\\
%&=& \partial_1(q_n) \wedge \partial_{n-1}(q_n)\\
%&\geq&
\geq \dfrac{a_{n-1}}{a_n}q_n.
$$ %\end{eqnarray*}
Hence
\begin{eqnarray*}
\dfrac{d(q_n)}{a_n} &\geq & \dfrac{a_{n-1}}{a_n^2}\dfrac{a_n+a_{n+1}}{2}\\
&=&\dfrac{a_{n-1}}{2}\left(  \dfrac{1}{a_n}+\dfrac{a_{n+1}}{a_n^2}\right)\to\infty(0+\infty)=\infty.
\end{eqnarray*}
Therefore there is no constant $C$ such that $q_n\in B(0,C a_n)$, for all $n\in N$. \qed

\subsection{On Hausdorff and Assouad dimensions of homogeneous distances}
In this section we provide two particular translation-invariant distances $d_1$ and $d_2$ on $\R$. Namely, the Hausdorff dimension of $(\R, d_1)$ is infinite  and the three dimensions (topological, Hausdorff, and Assouad)  of 
$(\R, d_2)$ are all different.
In particular, the metric space $(\R, d_1)$ shall be a non-doubling metric space. In \cite{LeDonne2}, it is showed that biLipschitz homogeneous geodesic planes are locally doubling.
The example $(\R, d_1)$ shows that the same conclusion does not hold for non-geodesic spaces, even when they are isometrically homogeneous and topological curves.

To define a translation-invariant
distance on $\R$ it is enough to define $\rho(r):=d(0,r)$, for $r\geq0$. However, not all functions $\rho:[0,\infty)\to[0,\infty)$ lead to distances.

First, we give a complete characterization of those functions $\rho$ that gives proper distances on $\R$ inducing the same topology. Second we present a sufficient condition for having such distances, which will be easier to check.

\begin{lemma}\label{characterization}
Let $X$ be a subgroup of the additive group $\R^n$.
Let $d:X\times X\to \R$ be a function.
Then $d$ is a proper $X$-invariant distance if and only if there exists a function $h:X\to \R$ such that
\begin{equation}d(x,y)=h(x-y),\qquad \forall x,y\in X,\label{d=hd_E}\end{equation}
satisfying the following properties:
\begin{description}
\item[\ref{characterization}.1] $h(0)\geq0$ and $h(x)=0$ if and only if $x= 0,$
\item[\ref{characterization}.2] $h(x)=h(-x), \qquad \forall x\in X,$
\item[\ref{characterization}.3] $h$ is sublinear, i.e.,
$$h(x+y)\leq h(x)+h(y),\qquad \forall x,y\in X,$$
\item[\ref{characterization}.4] $h$ is continuous at $0$ (and therefore, from the other properties, $h$ is uniformly continuous),
\item[\ref{characterization}.5] $$\inf_{|x|>\eps} h(x)>0, \qquad \forall \eps>0,$$
\item[\ref{characterization}.6] $$\sup_{h(x)<r} |x|< +\infty, \qquad \forall r>0.$$
\end{description}
\end{lemma}
\proof
Let $h$ be a function satisfying the six properties and such that  \eqref{d=hd_E} holds.
By  \eqref{d=hd_E}, the function $d$ is $X$-invariant.
The first two properties gives positive definiteness and symmetry, respectively.
The third property gives the triangle inequality. Indeed, for all $x,y,z\in X$,
$$d(x,y)=h(x-y)=h(x-z+z-y)\leq  h(x-z)+h(z-y)=
d(x,z)+d(z,y).$$	
Let us show that $d$ induces the usual topology, using the fourth and fifth property. Since $d$ is $X$-invariant it is enough to show that the neighborhoods of $0$ are the same. Let $\eps>0$, since $h$ is continuous at $0$, there exists $\delta>0$ such that
$$|x|<\delta \implies |h(x)|<\eps.$$
In other words,
$B_d(0,\eps)\supset B_{\rm usual}(0,\delta).$
On the other hand, since $\delta':=\inf_{|x|>\eps} h(x)>0$ is non-zero, we have that
$$|h(x)|<\delta' \implies |x|<\eps,$$
i.e., $B_d(0,\delta')\subset B_{\rm usual}(0,\eps).$
By the sixth property, for any $r>0$ we have $B_d(0,r)\subset B_{\R^n}(0,R),$
with $R:=\sup_{h(x)<r} |x|$. Therefore the metric $d$ is proper.

Viceversa, if $d$ is a  proper $X$-invariant distance, define $$h(x):=d(0,x), \qquad x\in X.$$
The equation  \eqref{d=hd_E} and the first three properties of $h$ are consequences of the fact that $d$ is a distance and that it is $X$-invariant.
Since $d$ induces the usual topology, we have that for any $\eps>0$ there exists $\delta>0$ such that
$B_d(0,\eps)\supset B_{\rm usual}(0,\delta)$, i.e.,
$h$ is continuous at $0$. In fact, $h$ is uniformly continuous since
$$ h(x)-h(y) \leq  h(x+y)\leq h(x)+h(y).$$
Again, since $d$ induces the usual topology, we have that for any $\eps>0$ there exists $\delta'>0$ such that  $B_d(0,\delta')\subset B_{\rm usual}(0,\eps)$, i.e.,
$\inf_{|x|>\eps} h(x)>\delta'>0.$
Finally, the last property comes from the properness of $d$, i.e., the fact that $\sup h^{-1}(0,r)<+\infty.$
\qed

\begin{corollary}\label{COROL-CHARACT}
Let $\rho:[0,\infty)\to[0,\infty)$ be a subadditive homeomorphism.
Then $d(x,y):=\rho(|x-y|)$, for $x,y\in\R$,  is a proper translation-invariant distance on $\R$.
\end{corollary}
\proof Define $h(x):=\rho(|x|),$ for $x\in\R$. All properties, except \ref{characterization}.3 of Lemma \ref{characterization} are obvious.
Regarding sublinearity, since $\rho$ is increasing and subadditive, we have
\begin{eqnarray*}
h(x+y)&=&\rho(|x+y|)\\
&\leq &\rho(|x|+|y|)\\
&\leq &\rho(|x|)+\rho(|y|)\\
&=&h(x)+h(y).
\end{eqnarray*}
\qed

\begin{remark}\label{COROL-CHARACT-rem}
Recall that, if a function $\rho:[0,\infty)\to[0,\infty)$ is twice differentiable with $\rho''\leq 0$, then $\rho$ is concave and therefore subadditive.
\end{remark}

For evaluating the Hausdorff measure, we will make use of Geometric Measure Theory arguments. The measure that we shall use is the Lebesgue measure, which we denote by $|A|$, given a  measurable set $A\subset \R$.
Next lemma tells us what is the measure of balls and it will repeatedly used later.

\begin{lemma}\label{measure ball}
Let $d$ be a translation-invariant metric on $\R$. Set $\rho(x):=d(0,x)$, for $x\geq0$.
Set $x_r:=\max\{x\;:\;\rho(x)=r\}$, for $r\geq0$.
Then we have
$$|B(p,r)| \leq 2 x_r.$$
If, moreover, the function $\rho:[0,\infty)\to[0,\infty)$ is strictly increasing, then, for all $p\in \R$ and $r\geq0$, we have
$$|B(p,r)| =2\rho^{-1}(r).$$
\end{lemma}
\proof
By translation invariance, we may assume $p=0$. 
Therefore, 
\begin{eqnarray*}
B(0,r)&=&\{x\in\R\;:\;d(0,x)<r\}\\
&=&\{x\in\R\;:\;\rho(|x|)<r\}\\
&\subseteq &(-x_r,x_r).
\end{eqnarray*}
If  $\rho $ is strictly increasing, we have
$$\{x\in\R\;:\;\rho(|x|)<r\}=
(-\rho^{-1}(r),\rho^{-1}(r)).$$
\qed

The following lemma is a collection of well-know facts from Geometric Measure Theory. The proofs of these results follow from Section 2.10.19 of \cite{Federer}.
\begin{lemma}\label{GMT}
Let $X$ be a metric space. Let $\mu$ be a Radon measure on $X$. Let $A\subseteq X$ be an open set of positive and finite measure.
\begin{description}
\item[\ref{GMT}.i] If there exist $\alpha, c>0$ such that $$\limsup_{r\to0}\dfrac{\mu (B(a,r))}{r^\alpha}=c, \qquad\forall a\in A, $$
then $\dim_H(A)=\alpha$.
\item[\ref{GMT}.ii]  If  $\alpha>0$ and $$\lim_{r\to0}\dfrac{\mu (B(a,r))}{r^\alpha}=0, \qquad\forall a\in A, $$
then the Hausdorff $\alpha$-dimensional measure $\mathcal H^\alpha(A)$ of $A$ is infinite.
\item[\ref{GMT}.iii]  If, for all $\alpha>0$, one has $$\lim_{r\to0}\dfrac{\mu (B(a,r))}{r^\alpha}=0, \qquad\forall a\in A, $$
then $\dim_H(A)=\infty$.
\end{description}
\end{lemma}

\begin{proposition}\label{Ex3}
Let $$\rho(x):=\min\left\{\sqrt{-\dfrac{1}{\log(x)}}\;,\; \left(\dfrac{2}{3}\right)^{ {3}/{2}}\left(\dfrac{e^{ {3}/{2}}}{2}x+1\right)\right\}.$$
Define $d(x,y):=\rho(|x-y|).$
Then $(\R,d)$ is a proper metric space that is isometrically homogeneous but whose Hausdorff dimension is infinite.
\end{proposition}
\proof By construction, the function $\rho$ is  a concave homeomorphism of $[0,\infty)$. Thus, by Corollary \ref{COROL-CHARACT} and Remark  \ref{COROL-CHARACT-rem}, the function $d$ is a proper translation-invariant distance (inducing the same topology).
By Lemma \ref{measure ball} and the fact that $\rho^{-1}(r)=e^{-1/r^2}$, for  $r$ sufficiently small, we have
$$\lim_{r\to0}\dfrac{ |B(a,r)|}{r^\alpha}=\lim_{r\to0}\dfrac{ e^{-1/r^2}}{r^\alpha}=0, \qquad \forall \alpha>0.$$ 
Thus, by \ref{GMT}.iii, the Hausdorff dimension of every ball is infinite. \qed

Regarding  the study of the Assouad dimension, we consider the following fact.
\begin{lemma}\label{lemma Assouad}
Let $(X,d,\mu)$ be a metric measure space.
Assume $$\mu(B(x,r))=\mu(B(x',r)), \qquad \forall x,x'\in X, \forall r>0.$$
Assume that there exists $\beta, c>0$ such that
$$\liminf_{r\to0}\dfrac{\mu (B(x,r))}{r^\beta}=c, \qquad\forall x\in X.$$
Then the Assouad dimension $\dim_A (X)$ of $X$ is $\geq \beta$.
\end{lemma}
\proof Let $b $ be a positive number with the property that there exists some $k$ such that
 every set of diameter $D$ can be covered by at most
 $k \eps^{-b}$ sets of diameter at most $\eps D$.   
 %, for all $x\in X, r>0, \eps>0$, the ball $B(x,r)$ can be covered with less than $k \eps^{-b}$  balls of radius $\eps r$.
Fix $x_0\in X$. Let $r_n\to 0$ such that
$$c/2\leq \dfrac{\mu (B(x_0,r_n))}{r_n^\beta}\leq 2c, \qquad \forall n\in \N.$$

Then, for $\eps:=r_m/r_n$, cover the ball $B(x_0,r_n)$, which has diameter less than $2r_n$,  with a collection of balls $B(x_1,  r_m), \ldots, B(x_N,  r_m)$, $N\in \N$. Since such balls have diameter at most $2 r_m=\eps(2r_n)$,  for the property of $b$ we can choose $N< k \eps^{-b}$.
Hence we have the bounds
\begin{eqnarray*}
\dfrac{c}{2} r_n^\beta &\leq & \mu (B(x_0,r_n))\\
&\leq & \sum_{j=1}^N \mu (B(x_j, {r_m} ))\\
&= & \sum_{j=1}^N \mu (B(x_0, r_m))\\
&\leq &  k \eps^{-b}  2 c r_m^\beta.
\end{eqnarray*}
Thus, for a suitable constant $c'$, we have 
$r_n^{\beta-b} \leq r_m^{\beta -b}$, for all $m$.
Since $r_m\to 0$ and $r_n>0$, we have $\beta -b<0$,
Hence
$\beta\leq\dim_A (X) $.
\qed

\begin{proposition}\label{Ex4}
Let  $\rho:[0,\infty)\to[0,\infty)$ be a concave homeomorphism 
such that $\sqrt[3]{x}\leq \rho(x)\leq \sqrt{x}$, for $x\in(0,1)$,
and for which there exist two sequences $x_n\to0$, $y_n\to0$ such that
$\rho(x_n) = \sqrt{x_n}$ and $\rho(y_n) = \sqrt[3]{y_n}$.
Define $d(x,y):=\rho(|x-y|).$
Then $(\R,d)$ is a proper metric space that is isometrically homogeneous for which
$$\dim_{top}(\R,d)\neq\dim_H(\R,d)\neq \dim_{Ass}(\R,d).$$
\end{proposition}
\proof
By Corollary \ref{COROL-CHARACT}, 
and Remark  \ref{COROL-CHARACT-rem}, the function $d$ is 
a proper translation-invariant distance. % (inducing the same topology).
Moreover,  the topology induced by $d$ is the usual one, hence we have $\dim_{top}(\R,d)=1$.
By Lemma \ref{measure ball}, 
$$\limsup_{r\to0}\dfrac{ |B(a,r)|}{r^2}=2.$$
Hence, by Lemma \ref{GMT}, $\dim_{H}(\R,d)=2$.
Finally, by Lemma \ref{lemma Assouad}, since
$$\liminf_{r\to0}\dfrac{ |B(x,r)|}{r^3}=2,$$ 
we have 
$\dim_{Ass}(\R,d)\geq 3$.
\qed

\begin{remark}
It is easy to give examples of   concave homeomorphisms $\rho:[0,\infty)\to[0,\infty)$  
such that $\rho(x)\geq x^2$, for $x>0$,
and for which there exist two sequences $x_n\to0$, $y_n\to0$ such that
$\rho(x_n) = \sqrt{x_n}$ and $\rho(y_n) = \sqrt[3]{y_n}$.
 %Indeed, ...
\end{remark}

\subsection{A non-Euclidean space with all dimensions $1$}
In this section we give an example of an isometrically homogeneous metric space that is not biLipschitz equivalent to the Euclidean line that has Assouad dimension one. In particular, also the topological, the Hausdorff and Nagata dimensions are $1$ as well.

By Theorem \ref{isom_homog_curves}, an isometrically homogeneous metric space is of the form  $X=(\R,d_\rho)$ with 
$d_\rho(s,t):=\rho(|s-t|),$
for some continuous function $\rho:[0,\infty)\to[0,\infty)$ with $\rho(0)=0$.
If the function $\rho$ is between two linear functions then the distance is biLipschitz homeomorphic to the Euclidean distance. Next lemma claims that the opposite implication holds true.

\begin{lemma}\label{linear rho}
Let $\rho:[0,\infty)\to[0,\infty)$ be such that, for $x,y\in\R$, $d_\rho(s,t):=\rho(|s-t|)$ is a distance function on $\R$, inducing the standard topology.
Assume that the metric space    $(\R,d_\rho)$ is biLipschitz homeomorphic to the Euclidean line.
Then there exists a constant $K$ such that
$$\dfrac{1}{K} x\leq \rho(x)\leq K x, \qquad \forall x>0.$$
\end{lemma}

\proof
Let $f: (\R,d_\rho)\to (\R,\norm{\cdot})$ be an $L$-biLipschitz homeomorphism, $L>1$.
We may assume that $f(0)=0$ and that $f(x)>0$, if $x>0$.
Hence, since $f$ is orientation preserving, if $a<b$, then $f(b)-f(a)>0$.

Thus, since $f$ is $L$-biLipschitz, we have that,  for $n,k\in \N$, 
\begin{eqnarray*}f\left(\dfrac{n}{k}\right)&=&\sum_{j=1}^n f\left(\dfrac{j}{k}\right)-f\left(\dfrac{j-1}{k}\right)\\
&\leq&\sum_{j=1}^n L \rho\left(\dfrac{1}{k}\right)=n L \rho\left(\dfrac{1}{k}\right).
\end{eqnarray*}
Similarly,
$$f\left(\dfrac{n}{k}\right)\geq n \dfrac{1}{L} \rho\left(\dfrac{1}{k}\right).
$$
In particular, 
for $n=k$, we have
$$\dfrac{1}{L}k \rho\left(\dfrac{1}{k}\right)\leq f\left(1\right)\leq L k  \rho\left(\dfrac{1}{k}\right).$$
Thus, for all $n,k\in \N$, we have
$$f\left(\dfrac{n}{k}\right)\leq n L \rho\left(\dfrac{1}{k}\right)\leq L^2 \dfrac{n}{k}f(1)$$
and
$$f\left(\dfrac{n}{k}\right)\geq n \dfrac{1}{L} \rho\left(\dfrac{1}{k}\right)\geq  \dfrac{1}{L^2}  \dfrac{n}{k}f(1).$$
Therefore,   for a suitable constant $C>0$, we have that, for all positive rational $x$,
$$C^{-1}x \leq f(x)\leq C x.$$
By continuity of $f$, the above equation holds for all $x\leq0$.

Using again that $f$ is $L$-biLipschitz, we get that
$$\rho(x)\leq L f(x)\leq CLx$$
and
$$\rho(x)\geq L^{-1} f(x)\geq (CL)^{-1}x.$$
Hence for $K:=CL$ the lemma is proved.
\qed

\begin{proposition} \label{dim1}
Let $\rho:[0,\infty)\to[0,\infty)$ such that 
$$
\rho^{-1}(x)= \left\{\begin{array}{lr}
-\int_0^x \dfrac{dt}{\log t} & \text{ for } \; x\in (0,1/2)\\
-\dfrac{1}{\log \half} (x-\half )-\int_0^{\half} \dfrac{dt}{\log t} &\text{ for } \; x\geq 1/2
.\end{array}\right.$$
Define $d(x,y):=\rho(|x-y|).$
Then $(\R,d)$ is a proper metric space that is isometrically homogeneous, has Assouad dimension $1$, Nagata dimension $1$, but is not locally biLipschitz homeomorphic to the Euclidean line.
\end{proposition}
\proof
The function $\rho^{-1}(x)$ vanishes at zero, is continuous, and increasing. Indeed,
$$
\dfrac{d}{dt}\rho^{-1}(x)= \left\{\begin{array}{lr}
-  \dfrac{1}{\log x} & \text{ for } \; x\in (0,1/2)\\
-\dfrac{1}{\log \half} &\text{ for } \; x\geq 1/2
.\end{array}\right.$$
Thus, $\rho^{-1}(x)$ is a $C^1$ homeomorphism of $[0,\infty)$.
We claim that the function $\rho$ is such that $\rho''\leq 0$. Indeed, we can show that
$\dfrac{d^2}{dt^2} \rho^{-1}(x)\geq 0$ .
In fact, for $ x\in (0,1/2)$,
$$\dfrac{d^2}{dt^2} \rho^{-1}(x)
=\dfrac{1/x}{\log^2x}
\geq 0.$$ 
Thus, by Corollary \ref{COROL-CHARACT} and Remark  \ref{COROL-CHARACT-rem}, the  pair $(\R,d)$ is a proper isometrically homogeneous metric space.
Since $\rho$ is increasing, the metric balls of $(\R,d)$ are connected. Hence, by Theorem \ref{dimension 1 intro},
the Nagata dimension of $(\R,d)$ is $1$.

Regarding the fact that such a metric space has Assouad dimension $1$, we need to show that, for all $\beta>1$, there exists a constant $c=c_\beta$ such that every set of diameter $D$, with $D>0$, can be covered by at most $c\eps^{-\beta}$ sets of diameter at most $\eps D$.
Thus we need to consider the set $[-\rho^{-1}(D),\rho^{-1}(D)]$ and cover it with translated of 
$[-\rho^{-1}(\eps D),\rho^{-1}(\eps D)]$.
Therefore, it is enough to show that the function
$$\eps \mapsto \dfrac{\eps^\beta \rho^{-1}(D)}{\rho^{-1}(\eps D)}
$$
is bounded in $\eps$, uniformly in $D$. It suffices to consider $D\in (0,1/2)$, since on large scale the distance is Euclidean.
In fact, the only issue might be as $\eps \to 0$. However, one sees that
$$ \dfrac{\eps^\beta \int_0^D \dfrac{dt}{\log t}}{\int_0^{\eps D} \dfrac{dt}{\log t}}=
\dfrac{\eps^\beta \int_0^D \dfrac{dt}{\log t}}{\int_0^{ D} \dfrac{\eps dt}{\log (\eps t)}}$$
is bounded in $\eps$, by showing that 
$ \dfrac{\eps^\beta /\log t }{\eps /\log (\eps t)}$ is bounded as $\eps \to 0$. Now, since $\beta-1>0$ and $t\in (0,1/2)$, we have
$$ \eps^{\beta-1} \dfrac{\log  (\eps t)}{ \log t}=\eps^{\beta-1} \dfrac{\log  \eps+\log  t}{ \log t}
\leq \eps^{\beta-1} \left(\dfrac{\log  \eps }{ \log (1/2)}+1 \right),$$
which tends to $0$, as $\eps \to 0$.

To show that such a metric space is not biLipschitz equivalent to the Euclidean line, one can either use Lemma 
\ref{linear rho} or prove that the $1$-dimensional Hausdorff measure is not locally finite. Indeed, the condition \ref{GMT}.ii is verified, because, by Lemma \ref{measure ball},
$$
\lim_{r\to 0} \dfrac{|B(r)|}{r}=
\lim_{r\to 0} \dfrac{2\rho^{-1}(r)}{r}=
\lim_{r\to 0}- \dfrac{2\int_0^r \dfrac{dt}{\log t} }{r}=
- \dfrac{2}{\log t}|_{t=0}=0.$$
\qed

%\begin{theorem} FALSE
%Let $X$
%be a metric space that is topologically equivalent to $\R$. Assume Isom$(X)\acts X$ transitively.
%Assume that the Hausdorff $1$-measure of balls of $X$ is finite.
%Then $X$ is locally biLipschitz equivalent to the Euclidean line. 
%\end{theorem}
%
%\proof
%By Theorem \ref{isom_homog_curves}, we may assume that  $X=(\R,d_\rho)$ with 
%$d_\rho(s,t):=\rho(|s-t|),$
%for some continuous function $\rho:[0,\infty)\to[0,\infty)$ with $\rho(0)=0$.
%We claim that there exists 

 \bibliography{general_bibliography}
\bibliographystyle{amsalpha}
 
 \vskip 1in

\parbox{3.5in}{Enrico Le Donne:\\
~\\
enrico.ledonne@math.ethz.ch\\ 
enrico.ledonne@msri.org
}

\end{document}